\documentclass[12pt]{article}
\usepackage{epsfig}
\usepackage{amsthm}
\usepackage{amsfonts}
\usepackage{amsmath}
\theoremstyle{plain}
\newtheorem{theorem}{Theorem}

\newcommand{\E}{\mathbb E}
\newcommand{\p}{\mathbb P}

\def\blackslug{\hbox{\hskip 1pt \vrule width 4pt height 8pt depth 1.5pt
\hskip 1pt}}
\def\QED{\quad\blackslug\lower 8.5pt\null\par}

\def\b0{\mathbf{0}}

\title{An Urn Model of Diaconis }
\author{D. Siegmund \\
Department of Statistics, Stanford University, Stanford, CA 94305 \\
dos@stat.stanford.edu 
\and B. Yakir 
\\Department of Statistics, Hebrew University, Jerusalem, Israel
\\msby@mscc.huji.ac.il }

\begin{document}

\maketitle



\begin{center}{SUMMARY}\end{center}
An urn model of Diaconis and some
generalizations are discussed.  A 
convergence theorem is proved that implies for Diaconis' model
that the empirical distribution of
balls in the urn converges with probability one
to the uniform distribution.

\section{Introduction}
Diaconis has formulated the
following simple urn model. 

\noindent
EXAMPLE 1.
Let $G$ be a finite group, generated
by $g_1, \ldots, g_r$. Initially, an urn contains $r$ balls, each
labeled by one of the generating elements. At times $n = r+1,r+2,
\cdots$ two balls are drawn with replacement from the urn. The
labels on these balls are multiplied to form a new group element.
A ball, bearing this element as its label, is then added to the
urn, increasing the number of balls in the urn by one. Let $X_{k}$
be the label indicator with respect to the $k$th ball (i.e., $X_k$
is a vector of length $|G|$, with a one placed in the coordinate
associated with the ball's label and zeros elsewhere.)  Let
$p_{g,n} = \sum_{k=1}^n I_{\{X_{g,k}=1\}}/n$ denote the relative
frequency of balls labeled $g$ when the total number of balls in
the urn is $n$.  As an application of Theorem 2 below, we
verify a conjecture of Diaconis, that
$p_{g,n} \rightarrow |G|^{-1}$, for all $g
\in G$, as $n \rightarrow \infty$ with probability one.

\noindent
EXAMPLE 2.  A special case of Example 1 occurs when the
balls are numbered either 0 or 1 and the group operation is
addition modulo 2.  Then $p_n,$ the
fraction of 1's in the urn after $n$ 
draws, converges to 1/2 with
probability one.   As a variation of this special case one
can draw $k \geq 2$ balls from the urn with replacement and add a
0 or a 1 according as the number of 1's drawn
is even or odd.  Again the fraction of balls
numbered 1 converges to 1/2 with probability one.  

\noindent
EXAMPLE 3.  For an example motivated by a classical model in
population genetics
(e.g. Ewens (1969)), we suppose that the population size in
a pure birth process at the $n$th generation is $k_n \geq n.$  The population
consists of three kinds of individuals corresponding to the 
three biallelic genotypes AA, Aa, and aa, which have relative
fitness (i.e., probability of reproduction ) of
$1-s, \; 1,\; 1-t$, respectively.  We assume
$s < 1,\; t < 1.$  In the most interesting special case
$0 < s < 1,\; 0< t < 1$, so the heterozygote Aa has the greatest
fitness.  Let $p_n$ denote the fraction of A alleles in the
population at the  $n$ generation.  Then under random mating
the relative proportions of AA, Aa and aa genotypes that
reproduce in the $n+1$st generation are
$p_n^2(1-s):2p_n(1-p_n):(1-p_n)^2(1-t).$
We assume that reproduction occurs independently of the population size
process.  Does the fraction $p_n$ converge and what is its limit?
In this example it is natural to assume that $k_n$ grows
exponentially, so that the number of balls added to the urn
in each generation is comparable to the number of balls
already in the urn.  One could also add this feature to
Examples 1 and 2.  

\section{Convergence to a fixed point}

Consider a finite set $G$. Let $G^*$ be the simplex of probability
distributions over $G$ and let $T : G^* \rightarrow G^*$ be a map
of the simplex into itself. The point $q \in G^*$ is a fixed point
of the transformation if $T(q) = q$. Below we investigate
almost-sure convergence of the stochastic sequence of empirical
distributions $\{p_n\}$, define by the recursion:
\[
p_{n+1}=\frac{k_n}{k_{n+1}}p_{n} + \frac{\sum_{i=k_{n}+1}^{k_{n + 1}}
X_{i}}{k_{n+1}} = \frac{k_0p_0 +\sum_{i=1}^{k_{n+1}}
X_i}{k_{n+1}},
\]
where $\{k_{n}\}$ is a monotone sequence of integer-valued random
variable (i.e.\ $k_{n+1} \geq k_n+1$, for all $n$), and $X_i$ is a
random vector that indicates an element from $G$. The integer
$k_0$ is positive and $p_0$ is a given initial distribution
vector. Consider the filtration ${\cal F}_n = \sigma\{X_1, \ldots,
X_{k_n},k_1, \cdots, k_n,  k_{n+1}\}$, 
for $n \geq 1$. We assume that, conditional
on ${\cal F}_n$,
\begin{equation} \label{eq:definition.X}
\sum_{i=k_{n}+1}^{k_{n+1}} X_{i} \sim \mbox{Multinomial}(T(p_n),
k_{n+1} - k_n),
\end{equation}
and identify sufficient conditions to ensure the convergence of
$p_n$ to a contracting (cf. assumption A1 below)
fixed point of the transformation $T$.

Our argument is a two-fold application of the almost
supermartingale convergence theorem of Robbins and Siegmund
(1971).  We begin with a statement of that theorem:

\begin{theorem} \label{th:1}
Let $Z_n, \xi_n, \zeta_n$ be non-negative random variables adapted
to the increasing sequence of $\sigma$-algebras ${\cal F}_n.$
Suppose that for each $n$,
$$\E(Z_{n+1} | {\cal F}_n) \leq Z_n + \xi_n - \zeta_n.$$
Then $\lim Z_n$ exists and is finite and $\sum \zeta_n < \infty$
almost surely on the event where $\sum \xi_n < \infty$.
\end{theorem}

Our main result relies on the following assumptions on the
transformation $T$, the sequence $\{k_n\}$ and the initial
distribution $p_0$:
\begin{description}
\item[A1:] The collection $Q = \{q_0, q_1, \cdots, q_J\}$ 
of fixed points of $T$ is
non-empty and finite and the fixed point $q_0$ is contracting, i.e., 
$\|T(p) - q_0\|  < \|p - q_0\|$, for all $p \in
G^* - Q$. 
The point $q_0$ may be in the interior of $G^*$, but all other
fixed points are on the boundary (i.e.\ their supports are proper subsets
of $G$).

\item[A2:] For all $j>0$ let $c_j$ be a vector with 0's
in those coordinates where $q_j$ has positive mass 
and 1's in those coordinates where $q_j$ has no mass.
Assume $c_j$ is not equal to the zero
vector (which  is equivalent to assuming that $q_j$ is on
the boundary of $G^*$).  Further assume that 
$\langle c_j, p_0\rangle > 0$ and for $p$ not orthogonal
to $c_j$,
$\liminf_{p
\rightarrow q_j} \langle c_j, T(p) \rangle/\langle c_j, p\rangle >
1$.

\item[A3:] The increasing sequence, $k_n$, of random  integers satisfies 
$k_{n+1}/k_n \leq C$, for all $n$ and for some
constant $C > 1$ such that $C-1 < \min \{||q_i - q_j||: i \neq j\}$.
\end{description}

\begin{theorem} \label{th:2}
Under the above assumptions, $p_n \rightarrow q_0$ with
probability one as $n \rightarrow \infty$.
\end{theorem}
\noindent {PROOF.}  The proof consists of applications of
Theorem 1 to (a) $Z_n = \|p_n - q_0\|^2$ and (b) $Z_n = 1/\langle
c_j, p_n\rangle $. Consider first case (a).  Let $\pi_{n+1} =
(k_{n+1}-k_{n})/k_{n+1}$ and define $\bar X_{n+1} =
\sum_{i=k_n+1}^{k_{n+1}} X_i/(k_{n+1}-k_n)$. 
Observe that
$p_{n+1} - q_0 = (1-\pi_{n+1})(p_n - q_0) + \pi_{n+1}(\bar{X}_{n+1} - q_0)$.
We take the the conditional
expectation given ${\cal F}_n$ of the squared norm of this identity and use 
the facts that
(i) $ \E(\bar X_{n+1} | {\cal F}_n) =  T(p_n)$
and (ii) the (conditional)  second moment of a random 
variable is the sum of its variance
and the square of its expectation.
Then by regrouping terms and using
the Cauchy-Schwarz inequality and conditions A1, A3
we see that
\begin{eqnarray*}
\lefteqn{\E(Z_{n+1} |{\cal F}_n) =Z_n -
2\pi_{n+1}(1-\pi_{n+1})\big[Z_n - \langle p_n-q_0, T(p_n) -
q_0\rangle\big]} \\ &&  +\pi_{n+1}^2 \big[\E\big(\|\bar
X_{n+1} - T(p_n)\|^2|{\cal F}_n\big) + \| T(p_n) - q_0\|^2 - Z_n]
\\ &\leq& Z_n - Z_n\frac{k_{n+1}-k_n}{C\cdot k_{n}}\big(1 - \frac{\|T(p_n) -
q_0\|}{\|p_n-q_0\|}\big)+\frac{k_{n+1}-k_n}{k_{n+1}^2}.
\end{eqnarray*}
Hence by A1 and Theorem 1, since
\[
\sum_{n=0}^\infty \frac{k_{n+1}-k_n}{k_{n+1}^2} \leq
\int_0^\infty\frac{dx}{x^2} < \infty,
\]
we see that with probability one, $\lim Z_n$ exists and is finite
and the negative terms of the process are summable. By the
non-negativity of the terms involved and by the fact that
\[
\sum_{n=0}^\infty \frac{k_{n+1}-k_n}{k_{n}} \geq
\int_{k_0}^\infty\frac{dx}{x} = \infty,
\]
we can conclude that either $Z_n \rightarrow 0$ or $\|T(p_n) -
q_0\|/\|p_n-q_0\| \longrightarrow_{n \rightarrow \infty} 1$.
However, only fixed points produce equality in the contraction
inequality. Consequently by A3, with probability one 
$p_n$ converges to  some $q_j \in Q$, the set of fixed points.

To eliminate the possibility that some $q_j$ with $j > 0$ is the limit,
we consider case (b): $Z_n = 1/\langle c_j, p_n\rangle $.
Indeed, we let $A_j = \{p_n \rightarrow q_j\}$ and show that $Z_n$
converges to a finite limit on $A_j$, which would
be a  contradiction unless
$\p(A_j) = 0$. This will complete the proof of
the theorem since $p_n$ must converge to a fixed point.

We turn to proving the convergence of $\{Z_n\}$ on $A_j$. Define
$\tilde S_{n+1}= \langle c_j, \sum_{i=k_n+1}^{k_{n+1}} X_i
\rangle$, $\tilde p_n = \langle c_j, p_n \rangle$ and $\tilde
T(p_n) = \langle c_j, T(p_n) \rangle$. Note that $\tilde
p_{n+1} = [k_n \tilde p_{n} +(k_{n+1}-k_n)\tilde S_{n+1}]/k_{n+1}$.
Conditional on ${\cal F}_n$,
$\tilde S_{n+1}$ is the sum of a subset of the coordinates
of a multinomial vector and hence
is distributed as Binomial(
$k_{n+1}- k_n, \tilde T(p_n)$). Now
\[
\E[Z_{n+1}| {\cal F}_n]= \E\bigg[\frac{k_{n+1}}{k_n \tilde
p_n+\tilde S_{n+1}}\bigg | {\cal F}_n\bigg] = \sum_{s=0}^{k_{n+1}
- k_n} \frac{k_{n+1}}{k_n\tilde p_n+s} \p(\tilde S_{n+1}=s | {\cal
F}_n).
\]
The relations $\p(\tilde S_{n+1}=0|{\cal F}_n) = 1 -
\sum_{s=1}^{k_{n+1}-k_n} \p(\tilde S_{n+1}=s|{\cal F}_n)$ and
$1/(k_n \tilde p_n+s) - 1/(k_n\tilde p_n) = -s/(k_n \tilde
p_n+s)\cdot1/(k_n \tilde p_n)$ produces
\begin{equation}\label{eq:sub.martingale}
= Z_n + \frac{k_{n+1}-k_n}{k_n \tilde p_n}\bigg[1
-\frac{k_{n+1}}{k_{n+1}-k_n}\sum_{s=1}^{k_{n+1}-k_n} \frac{s \cdot
\p(\tilde S_{n+1}=s |{\cal F}_n)}{k_n \tilde p_n+s} \bigg].
\end{equation}
We will proceed by showing that on the event $\{\tilde p_n
\rightarrow 0\} \supset A_j$ the term in the square brackets is
eventually strictly negative. Therefore, the positive part is
summable, and Theorem~\ref{th:1} can be used in order to conclude
that $\lim Z_n$ exists and is finite.

We analyze separately the cases: (i) $\E( \tilde S_{n+1} |{\cal
F}_n) < \epsilon$, (ii) $\epsilon \leq \E( \tilde S_{n+1} |{\cal
F}_n) \leq M$, and (iii) $\E( \tilde S_{n+1} |{\cal F}_n)
> M$, for some prespecified $0 < \epsilon < M < \infty$ 
to be determined later.

Consider case (i). By the monotonicity of the function $x/(a+x)$
we obtain the inequality
\[
\bigg [\cdots \bigg] \leq \bigg[1
-\frac{k_{n+1}}{k_{n+1}-k_n}\frac{\p(\tilde S_{n+1}\geq 1 | {\cal
F}_n)}{k_n \tilde p_{n}+1} \bigg].
\]
Now, $\p(\tilde S_{n+1}\geq 1|{\cal F}_n) = 1- (1- \tilde T(
p_n))^{k_{n+1}-k_n} \geq (k_{n+1}-k_n) \tilde T(p_n)(1-
\epsilon/2)$, which leads to the inequality
\[
\leq \bigg[1 -(1- \epsilon/2)\frac{k_n T(\tilde p_n)}{k_n \tilde
p_n+1}\bigg].
\]
If $k_n \tilde p_n \rightarrow \infty$, then assumption A2 will
produce a negative limit provided that $\epsilon$ is small enough.

To prove that $k_n \tilde p_n \rightarrow \infty$, it is
sufficient to prove that $\sum_{n=0}^\infty I_{\{\tilde S_{n+1}
\geq 1\}}$ is almost surely infinite. Equivalently, it is enough
to show
\[
\sum_{n=n_0}^\infty \p(\tilde S_{n+1} \geq 1 |{\cal F}_n) \geq
\sum_{n=n_0}^\infty (k_{n+1}-k_n)\tilde p_n(1- \epsilon/2) =
\infty,
\]
for an appropriate $n_0$. However, $\tilde p_n \geq \langle c_j,
p_0\rangle/k_n$, and the statement follows from the fact that
$\{(k_{n+1} - k_n)/k_n\}$ has an infinite sum.

Next consider case (ii). Since $\tilde T(p_n) \rightarrow
0$, we must have that $k_{n+1}-k_n \rightarrow \infty$ and thus
$\tilde S_{n+1}$ behaves in distribution like a Poisson random
variable (conditional on ${\cal F}_n$). This time we use the
inequality
\[
\bigg [\cdots \bigg]  \leq \bigg[1 -\frac{1}{(k_{n+1}-k_n)\tilde
p_n}\E \bigg(\frac{\tilde S_{n+1}}{1+\tilde S_{n+1}/k_n\tilde p_n}
\bigg | {\cal F}_n\bigg)\bigg].
\]
Case (ii) implies a lower bound on the term $(k_{n+1}-k_n)\tilde
p_n$ and a stochastic upper bound on the random variable $\tilde
S_{n+1}$. It follows that the conditional expectation $\sim \E (\tilde
S_{n+1} | {\cal F}_n)= (k_{n+1}-k_n)\tilde T(p_n)$, which
produces a negative value in the square brackets, by A2.

Finally, consider case (iii). By monotonicity one gets
that
\[
\frac{s}{a+s}\geq \frac{y \cdot I_{\{s \geq y\}}}{a+y}
\]
and, upon selecting $y = (1-\epsilon_1)\E(\tilde S_{n+1} |{\cal
F}_n)$, the inequality
\[
\bigg [\cdots \bigg] \leq \bigg[1 -\frac{k_{n+1}\p\big(\tilde
S_{n+1} \geq (1-\epsilon_1)\E(\tilde S_{n+1} |{\cal F}_n)|{\cal
F}_n\big)}{k_n[\tilde p_n/(1-\epsilon_1)\tilde T(
p_n)]+(k_{n+1}-k_n)} \bigg].
\]
Chernoff's inequality leads to the upper bound
\[
\bigg[1 -\frac{k_{n+1}}{k_n[\tilde p_n/(1-\epsilon_1)\tilde
T(\tilde p_n)]+(k_{n+1}-k_n)}  (1- e^{-\frac{\epsilon_1^2
M}{2}})\bigg].
\]
Selection of a large enough $M$ and a small enough $\epsilon_1$
will lead to a negative limit, provided that
$(k_{n+1}-k_n)/k_n$ is bounded. This last condition is assured
by  assumption A3.

\section{Applications}

\noindent
EXAMPLE 1.
In the urn model of Diaconis the transformation takes the form:
\[
\big(T(p)\big)_g = \sum_{h \in G} p_{g\cdot h^{-1}} p_{h}, \quad
\mbox{for $g \in G$}.
\]
Any uniform distribution over a subgroup is a fixed point of this
transformation. Conversely, any fixed point is a uniform
distribution over a subgroup. The last statement follows from the
fact that the support of a fixed point is a subgroup since the
support is closed under group operations and the group is finite.
Moreover, by the definition of a fixed point, the probability of
each element in the support must be equal to the maximum of all
probabilities unless a contradiction is to occur. The collection
of uniform distributions over subgroups is finite.

Denote by $q_0$ the uniform distribution over the entire group.
Viewing $\big(\sum_{h \in G} p_{g\cdot h^{-1}} p_{h}\big)^2$
as the square of the expectation of the random variable
taking on the value $p_h$ with probability $p_{g\cdot h^{-1}}$,
we obtain from the 
Cauchy-Schwarz inequality that $\sum_{g \in G}
\big(\sum_{h \in G} p_{g\cdot h^{-1}} p_{h}\big)^2\leq \sum_{g \in
G}p_{g}^2$, with strict
inequality unless $p_h$ is constant on its support.  From this 
and direct computations, we see that
$T$ is contracting, so condition A1 is met.

Let $G_j$ be a proper sub-group of $G$. Observe that $\langle c_j,
p\rangle$ assigns a probability to $G \setminus G_j$. A product of
two group elements, one belonging to $G_j$ and the other not
belonging, produces a group element not belonging to $G_j$. It
follows that
\[
\langle c_j, T(p)\rangle \geq 2 \langle c_j, p\rangle(1-\langle
c_j, p\rangle).
\]
If $p_0$ assigns positive probabilities to generators of $G$ then
$\langle c_j, p_0\rangle > 0$ and condition A2 is fulfilled.

\noindent
EXAMPLE 2.
From the elementary fact that when a coin is
tossed $k$ times, the probability of an
odd number of heads is $[1-(1-2p)^k]/2$, one can
verify the conditions of the theorem, to show that
$p_n \rightarrow 1/2$ with probability one.  It is perhaps
interesting to note that when $k$ is even the transformation
$T(p)$ is concave;  when $k$ is odd, it is concave to the left
of 1/2 and convex to the right of 1/2.

\noindent
EXAMPLE 3.  From the assumption of random mating
it follows that $T(p) = p(1-ps)/[1-p^2s-(1-p)^2t]$, from
which it easily follows that 0 and 1 are fixed points 
of $T$.  If $s$ and $t$ are both positive or both 
negative, then $q^* = t/(s+t)$ is also a fixed point;
otherwise 0 and 1 are the only fixed points.
It is straightforward to show that when $s$ and $t$ are both positive,
the interior point $t/(s+t)$ is attracting, so $p_n \rightarrow
t/(s+t)$ with probability one.  (Like Example 2,
$T$ is concave to the left of $q^*$ and 
convex to the right.)  When $s$ is nonpositive and $t$
is positive, the fixed point at 1 is attracting, and conversly
in the case when $s$ is positive and $t$ nonpositive.  
If $s = t = 0$, every point in 
[0,1] is a fixed point, the sequence $p_n$ is a martingale,
which converges with probability one to a random limit.
In the
case when both $s$ and $t$ are negative, the fixed point at
$t/(s+t)$ is not attracting.
It seems intuitively clear that
$p_n$ must converge to 0 or 1, but this does not
seem to follow from Theorem 2 without an additional
argument.

\medskip\noindent
ACKNOWLEDGMENT.
The authors would like to thank Persi Diaconis for suggesting
this problem in the first place, Steve Lalley for
suggesting we find a more general formulation of our
original result, which
dealt only with Example 1, and a referee whose careful reading
has prevented us from making at least one egregious error
in the formulation of Theorem 2.
This research has been partially supported by the National
Science Foundation and by the U.S.-Israel Binational Science
Foundation.

\begin{center} {REFERENCES} \end{center}
\begin{description}

\item Duflo, M. (1993).   {\sl Random Iterative Models},
Springer-Verlag, New York.

\item Ewens, W. (1969).  {\sl Population Genetics}, 
Methuen, London.

\item Robbins, H. and Siegmund, D. (1971).  A convergence
theorem for non-negative
almost supermartingales
and some applications, in {\sl Optimizing Methods in
Statistics}, Academic Press, New York, pp.  233--257.

\end{description}

\end{document}